# Results for Global Attractivity of Interior Equilibrium Points for Lotka-Volterra Systems

**Iasson Karafyllis**

*Dept. of Mathematics, National Technical University of Athens, Zografou Campus, 15780, Athens, Greece.
emails: iasonkar@central.ntua.gr , iasonkaraf@gmail.com

## Abstract

This paper provides global attractivity results for the interior equilibrium point of a general Lotka-Volterra system with no restriction on the dimension of the system and with no special structure or properties of the interaction matrix. The main result contains as special cases all known general results, including the Volterra-Lyapunov theorem and the recently proposed eigenvector conditions. Moreover, global attractivity of the interior equilibrium point is shown for a three-dimensional example, where none of the existing general results can be applied.

## 1. Introduction

Lotka-Volterra models are popular because of their simplicity and their efficiency to describe qualitatively various phenomena in ecology. Many researchers have studied dynamical properties of Lotka-Volterra models, including boundedness of solutions, persistence, permanence, existence and stability of equilibria, limit cycles and heteroclinic orbits (see for instance [4, 15]).

Interior equilibrium points for Lotka-Volterra models, i.e., equilibrium points where all populations coexist, are important (both in mathematics and ecology). The stability properties of interior equilibrium points have been studied extensively in the literature. Special results have been obtained for competitive Lotka-Volterra models, where the existence of a carrying simplex, i.e., an attracting invariant manifold that determines the dynamical behavior, can be proved (see [3, 10]). The split Lyapunov function method was used in [7, 16, 22] in order to prove global attractivity of the interior equilibrium point while in [5, 6] the notion of the ultimate contracting cells was exploited. The three-dimensional competitive case has also attracted the attention of researchers: in [11] the authors provided criteria for global stability based on the geometric approach using compound matrices, while in [20] the authors gave criteria for global stability and the non-existence of limit cycles.

The global stability properties of the interior equilibrium point have been studied for other special classes of Lotka-Volterra systems (including predator-prey structures): see [2, 17, 18, 19].

For general Lotka-Volterra models, i.e., models that do not assume a special structure or special properties for the interaction matrix, there are two major results in the literature that guarantee global attractivity of an interior equilibrium point:

**(A)** The Volterra-Lyapunov conditions (see [4]). The classical Volterra-Lyapunov theorem has been extended in [12, 13] by using LaSalle's theorem.

**(B)** The eigenvector conditions (see [1]). More specifically, the authors in [1] extended the computational theorem in [22] to the general (not necessarily competitive) case.

There are also global attractivity results in [21] that deal specifically with the three-dimensional case. The case studied in [21] includes Lotka-Volterra systems which are not necessarily competitive.

The present work provides global attractivity results for the interior equilibrium point of a general Lotka-Volterra system with no restriction on the dimension of the system. The main result (Theorem 1) contains as special cases all known general results: the Volterra-Lyapunov theorem, the eigenvector conditions in [1] and the global attractivity results in [21]. Moreover, we study a three-dimensional example (suggested in [1] as an open problem), where none of the aforementioned results can be applied and where global attractivity of the interior equilibrium point can be proved by using the results of the present work.

The methodology followed in the present work is the construction of a family of functions that generalizes the well-known class of Volterra-Lyapunov functions. The functions are used in conjunction with LaSalle's theorem on specific trapping regions.

The structure of the paper is as follows. Section 2 clarifies the notions that are used in the paper and provides the main results of the paper. Section 3 of the present work is devoted to the investigation of the relation of the obtained results with existing results in the literature. The proof of the main result (Theorem 1) is provided in Section 4. Finally, some concluding remarks and directions for future work are given in Section 5.

**Notation.** Throughout this paper we use the following notation.

- Consider the dynamical system $\dot{x} = f(x), x \in D$, where $D \subseteq \mathbb{R}^n$ is an open set and $f : D \to \mathbb{R}^n$ is a locally Lipschitz mapping. For a continuously differentiable function $V : D \to \mathbb{R}$ we define $\dot{V}(x) = \nabla V(x) f(x)$ for all $x \in D$. Let $x_0 \in D$ be such that the solution $x(t)$ of $\dot{x} = f(x)$ with initial condition $x(0) = x_0$ exists for all $t \geq 0$. We define $\omega(x_0) \subseteq \mathbb{R}^n$, the $\omega-$ limit set of $x_0 \in D$ using the formula:

$$\omega(x_0) = \left\{ y \in \mathbb{R}^n : \exists \{t_k \geq 0 : k = 1, 2, ...\} \text{ with } \lim_{k \to +\infty} (t_k) = +\infty \text{ and } \lim_{k \to +\infty} \left( x(t_k) \right) = y \right\}$$

- For a non-empty set $A \subseteq \mathbb{R}^n$ we denote its closure by $\bar{A} \subseteq \mathbb{R}^n$, i.e., $\bar{A} = \left\{ x \in \mathbb{R}^n : dist(x, A) = 0 \right\}$, where $dist(x, A) = \inf \left\{ |y - x| : y \in A \right\}$. The interior of a set $A \subseteq \mathbb{R}^n$ is denoted by $\mathrm{int}(A)$.

- For a vector $x \in \mathbb{R}^n$, $|x|$ denotes its Euclidean norm and $x^T$ denotes its transpose. If $B \subseteq \mathbb{R}$, the set of $n \times n$ matrices with entries in $B \subseteq \mathbb{R}$ is denoted by $B^{n \times n}$. For a matrix $C \in \mathbb{R}^{n \times n}$ we denote its transpose by $C^T$. We define $C^S = \frac{1}{2}\left(C + C^T\right)$ for all $C \in \mathbb{R}^{n \times n}$. A symmetric matrix $C \in \mathbb{R}^{n \times n}$ is called positive semidefinite if $x^T C x \geq 0$ for all $x \in \mathbb{R}^n$. A symmetric matrix $C \in \mathbb{R}^{n \times n}$ is called positive definite if $x^T C x > 0$ for all $x \in \mathbb{R}^n \setminus \{0\}$. A symmetric matrix $C \in \mathbb{R}^{n \times n}$ is called negative (semi)definite if $(-C)$ is positive (semi)definite. For a vector $x \in \mathbb{R}^n$, $diag(x)$ denotes the diagonal $n \times n$ matrix with the entries of $x$ on its main diagonal.

## 2. Main Results

In this section we present all notions for Lotka-Volterra systems that are used in this paper and our main result. All notions for Lotka-Volterra systems are taken from [4].

Lotka-Volterra systems are dynamical system of the form

$$\dot{y} = diag(y)\left(r + By\right) \tag{1}$$

where $B \in \mathbb{R}^{n \times n}$ is a constant real matrix, $r \in \mathbb{R}^n$ is a constant vector and $y \in \mathbb{R}^n_+$ is the state. It is well-known that $\mathbb{R}^n_+$, $\mathrm{int}\left(\mathbb{R}^n_+\right)$ and $\partial \mathbb{R}^n_+ = \mathbb{R}^n_+ \setminus \mathrm{int}\left(\mathbb{R}^n_+\right)$ are invariant sets for Lotka-Volterra systems.

A Lotka-Volterra system of the form (1) is called:

- uniformly bounded if there exists a constant $M > 0$ such that all solutions $y(t)$ of (1) exist for all $t \geq 0$ and satisfy the following inequality for all $i = 1,...,n$:

$$\limsup_{t \to +\infty}\left(y_i(t)\right) \leq M \tag{2}$$

- strongly persistent if all solutions $y(t)$ of (1) with $y(0) \in \mathrm{int}\left(\mathbb{R}^n_+\right)$ exist for all $t \geq 0$ and satisfy the following inequality for all $i = 1,...,n$:

$$\liminf_{t \to +\infty}\left(y_i(t)\right) > 0 \tag{3}$$

- permanent if (1) is uniformly bounded and there exists a constant $\delta > 0$ such that all solutions $y(t)$ of (1) with $y(0) \in \mathrm{int}\left(\mathbb{R}^n_+\right)$ satisfy the following inequality for all $i = 1,...,n$:

$$\liminf_{t \to +\infty}\left(y_i(t)\right) \geq \delta \tag{4}$$

A point $y^* \in \mathrm{int}\left(\mathbb{R}^n_+\right)$ is called an interior equilibrium point for (1) if $r + By^* = 0$. An interior equilibrium point $y^*$ is globally attractive for (1) if $\omega\left(y\right) = \left\{ y^* \right\}$ for all $y \in \mathrm{int}\left(\mathbb{R}^n_+\right)$.

When there is an interior equilibrium point $y^*$, the transformation $y = diag\left(y^*\right)x$ brings the Lotka-Volterra system (1) to the form:

$$\begin{aligned} &\dot{x} = diag(x)A\left(x - 1_n\right) \\ &x \in \mathbb{R}^n_+ \end{aligned} \tag{5}$$

where $1_n = \left(1,...,1\right)^T$ and $A = Bdiag\left(y^*\right)$. We next focus on system (5) that always has the interior equilibrium point $x^* = 1_n$.

In the literature there are two major results that guarantee global attractivity of an interior equilibrium point $x^* = 1_n$ for the general case of system (5), i.e., without assuming any special structure and properties for the matrix $A$.

**(A)** Volterra-Lyapunov conditions (see [4]). If there exists $h \in \mathrm{int}\left(\mathbb{R}_+^n\right)$ such that $\left(diag(h)A\right)^S$ is negative definite then the interior equilibrium point $x^* = 1_n$ is globally attractive for (5).

**(B)** Eigenvector conditions (see [1]). If
**(i)** system (5) is permanent,
**(ii)** there exists $\alpha = \left(\alpha_1,...,\alpha_n\right)^T \in \mathbb{R}^n$ and $\lambda < 0$ with $\alpha_i \neq 0$ for $i = 1,...,n$ and $\alpha^T A = \lambda \alpha^T$, and
**(iii)** for all $y \in \mathbb{R}^n$ with $y \neq 0$ and $\alpha^T y = 0$, it holds that $y^T diag(\alpha) A y < 0$
then the interior equilibrium point $x^* = 1_n$ is globally attractive for (5).
The eigenvector conditions in [1] extend the global attractivity results in [22] that deal with competitive systems, i.e., systems of the form (5) with $A \in \left(-\infty, 0\right)^{n \times n}$.

There are also additional global attractivity results in [21] that deal specifically with the case $n = 3$.

Using LaSalle's theorem and Lyapunov arguments we can prove the following theorem that provides sufficient conditions for convergence of a solution of (5) to the interior equilibrium point $x^* = 1_n$ for (5).

**Theorem 1:** *Suppose that there exist constants $\kappa, \mu > 0$, $g, \delta \geq 0$, $b \in \mathbb{R}$ and vectors $p, q, k \in \mathbb{R}^n$, $\beta \in \mathbb{R}_+^n$ such that the matrices $R = \left(\left(diag\left(k\right) + bkk^T\right)A\right)^S$ and $Q = \left(diag\left(q\right)A - \delta pp^T\right)^S$ are negative semidefinite. Moreover, suppose that*

$$\mu \geq \beta^T 1_n - \kappa g \tag{6}$$

$$\left(diag\left(p\right)A + \left(\beta - gp\right)p^T\right)^S = 0 \tag{7}$$

$$A^T p = -\left(\mu - g\kappa\right)p - \kappa\beta \tag{8}$$

*Let $F \subseteq \mathbb{R}_+^n$ be a closed set with the property that for every $x_0 \in \mathrm{int}\left(\mathbb{R}_+^n\right)$ for which the solution $x(t)$ of (5) with initial condition $x(0) = x_0$ exists for all $t \geq 0$, there exists $\tau > 0$ such that $x(t) \in F$ for all $t \geq \tau$. Finally, suppose that $\left\{1_n\right\} = \mathrm{int}\left(\mathbb{R}_+^n\right) \cap J$, where $J \subseteq F$ is the largest invariant set in the set*

$$C := \left\{ x \in F : \left(x - 1_n\right)^T R\left(x - 1_n\right) = \left(x - 1_n\right)^T Q\left(x - 1_n\right) = p^T\left(x - 1_n\right) = 0 \right\} \tag{9}$$

*Then the following statements hold:*

**(a)** *for every* $x_0 \in \operatorname{int}\left(\mathbb{R}_+^n\right)$ *for which the solution* $x(t)$ *of (5) with initial condition* $x(0)=x_0$ *is bounded and satisfies* $\liminf_{t\to+\infty}\left(x_i(t)\right)>0$ *for* $i=1,...,n$, *it holds that* $\omega\left(x_0\right)=\{1_n\}$,

**(b)** *if* $b=0$, $(-p)\in\mathbb{R}_+^n$ *and there exists* $c^*>0$ *such that* $q_i+k_i>c^* p_i$ *for* $i=1,...,n$ *then for every* $x_0 \in \operatorname{int}\left(\mathbb{R}_+^n\right)$ *for which the solution* $x(t)$ *of (5) with initial condition* $x(0)=x_0$ *is bounded it holds that* $\omega\left(x_0\right)=\{1_n\}$,

**(c)** *if* $k\in\mathbb{R}_+^n$, $(-p)\in\mathbb{R}_+^n$, $b<0$ *with* $bk^T 1_n>-1$ *and there exists* $c^*>0$ *such that* $q_i+k_i>c^* p_i$, $q_i\geq c^* p_i$ *for* $i=1,...,n$ *then for every* $x_0 \in \operatorname{int}\left(\mathbb{R}_+^n\right)$ *for which the solution* $x(t)$ *of (5) with initial condition* $x(0)=x_0$ *is bounded it holds that* $\omega\left(x_0\right)=\{1_n\}$.

The proof of Theorem 1 is provided in Section 4 and is based on the exploitation of the properties of the following parameterized family of functions with parameter $c>0$ which is defined for all $x\in\operatorname{int}\left(\mathbb{R}_+^n\right)$ with $p^T\left(x-1_n\right)+\kappa>0$:

$$U_c(x)=\sum_{i=1}^{n}\left(q_i+k_i-cp_i\right)f\left(x_i\right)+c\kappa f\left(\frac{p^T\left(x-1_n\right)+\kappa}{\kappa}\right), \text{ when } b=0$$

and

$$\begin{aligned} U_c(x)=&\sum_{i=1}^{n}\left(q_i-cp_i\right)f\left(x_i\right)+c\kappa f\left(\frac{p^T\left(x-1_n\right)+\kappa}{\kappa}\right)\\ &+\frac{1}{b}+\left(k^T\left(x-1_n\right)-\frac{1}{b}\right)\exp\left(b\sum_{i=1}^{n}k_i\ln\left(x_i\right)\right)\end{aligned}, \text{ when } b\neq 0$$

where $f(s)=s-1-\ln(s)$ for $s>0$. It should be noticed that when $b=0$ and $p=0$ the above family gives the well-known Volterra-Lyapunov functions $V(x)=\sum_{i=1}^{n}\left(q_i+k_i\right)f\left(x_i\right)$ which are utilized in the proof of Volterra-Lyapunov stability (see [4]).

Theorem 1 allows us to obtain the following corollaries which guarantee global attractivity of the interior equilibrium point $x^*=1_n$ for (5).

**Corollary 1:** *If*
**(i)** *the assumptions of Theorem 1 hold, and*
**(ii)** *system (5) is strongly persistent and for every* $x_0 \in \operatorname{int}\left(\mathbb{R}_+^n\right)$ *the solution* $x(t)$ *of (5) with initial condition* $x(0)=x_0$ *is bounded,*
*then the interior equilibrium point* $x^*=1_n$ *is globally attractive for (5).*

**Corollary 2:** *If*
**(i)** *the assumptions of Theorem 1 hold,*
**(ii)** $b=0$, $(-p)\in\mathbb{R}^n_+$ *and there exists* $c^*>0$ *such that* $q_i+k_i>c^* p_i$ *for* $i=1,...,n$, *and*
**(iii)** *for every* $x_0\in\text{int}\left(\mathbb{R}^n_+\right)$ *the solution* $x(t)$ *of (5) with initial condition* $x(0)=x_0$ *is bounded,*
*then the interior equilibrium point* $x^*=1_n$ *is globally attractive for (5).*

**Corollary 3:** *If*
**(i)** *the assumptions of Theorem 1 hold,*
**(ii)** $k\in\mathbb{R}^n_+$, $(-p)\in\mathbb{R}^n_+$, $b<0$ *with* $bk^T 1_n>-1$ *and there exists* $c^*>0$ *such that* $q_i+k_i>c^* p_i$, $q_i\geq c^* p_i$ *for* $i=1,...,n$, *and*
**(iii)** *for every* $x_0\in\text{int}\left(\mathbb{R}^n_+\right)$ *the solution* $x(t)$ *of (5) with initial condition* $x(0)=x_0$ *is bounded,*
*then the interior equilibrium point* $x^*=1_n$ *is globally attractive for (5).*

A detailed discussion of the significance of Theorem 1 and its relation with existing results in the literature is provided in the following section.

## 3. Examples and Comparison with Existing Results

In order to show the significance of Theorem 1, we start with an example of a Lotka-Volterra system for which *none of the existing general results in the literature* can be used for the proof of global attractivity of the unique interior equilibrium point. On the other hand, Theorem 1 can be used in a straightforward way and can allow us to conclude that the interior equilibrium point is globally attractive. The example was suggested in [1] as an open problem.

**Example 1:** System (1) with $n=3$ and

$$r=\begin{bmatrix}1\\3\\19/6\end{bmatrix},\ B=\begin{bmatrix}-1&0&9\\-1&-1&0\\-1&-1&-1\end{bmatrix}$$

was studied in [1]. The interior equilibrium point is $y^*=\begin{bmatrix}\frac{5}{2}&\frac{1}{2}&\frac{1}{6}\end{bmatrix}^T$ and the system takes the form (5) with

$$A=\frac{1}{2}\begin{bmatrix}-5&0&3\\-5&-1&0\\-5&-1&-1/3\end{bmatrix}\tag{10}$$

It was shown in [1] that the system is permanent but there is no $h\in\text{int}\left(\mathbb{R}^3_+\right)$ such that $\left(diag(h)A\right)^S$ is negative definite: the Volterra-Lyapunov stability theorem is not applicable. Moreover, the eigenvector conditions proposed in [1] are not applicable either: although there exists $\alpha = (\alpha_1, \alpha_2, \alpha_3)^T \in \mathbb{R}^3$ and $\lambda < 0$ with $\alpha_i \neq 0$ for $i = 1,2,3$ and $\alpha^T A = \lambda \alpha^T$, it is not true that $y^T diag(\alpha) Ay < 0$ for all $y \in \mathbb{R}^3$ with $y \neq 0$ and $\alpha^T y = 0$. Furthermore, it was mentioned in [1] that there is numerical evidence that the interior equilibrium point $x^* = 1_3$ is globally attractive.

Therefore, *none of the existing general results in the literature* can be used for the proof of global attractivity of the unique interior equilibrium point $x^* = 1_3$ for (5) with $A$ given by (10).

It is also clear that the results in [21] are also not applicable to this example since matrix $A$ does not have the special form (up to a renumbering of the states)

$$A = -\begin{bmatrix} \frac{y_1^*}{\lambda_1} & \frac{y_2^*}{\lambda_1} & \frac{\lambda_1 + \delta}{\lambda_1 \delta} y_3^* \\ \frac{y_1^*}{\lambda_2} & \frac{y_2^*}{\lambda_2} & \frac{y_3^*}{\lambda_2} \\ \frac{\delta - \lambda_3}{\lambda_3 \delta} y_1^* & \frac{y_2^*}{\lambda_3} & \frac{y_3^*}{\lambda_3} \end{bmatrix} \tag{11}$$

with

$$\begin{aligned} & 0 < \lambda_1 < \lambda_2 < \lambda_3 \\ & 0 < \delta < \lambda_1 \lambda_3 \frac{(1 - \lambda_2)}{\lambda_2 (\lambda_3 - \lambda_1)} \end{aligned} \tag{12}$$

and

$$\begin{aligned} & y_1^* = \frac{\delta(\lambda_3 - \lambda_2)}{\lambda_3} > 0 \\ & y_2^* = 1 - \lambda_2 - \delta \lambda_2 \frac{(\lambda_3 - \lambda_1)}{\lambda_1 \lambda_3} > 0 \\ & y_3^* = \frac{\delta(\lambda_2 - \lambda_1)}{\lambda_1} > 0 \end{aligned} \tag{13}$$

for which the results in [21] are applicable. Indeed, notice that the matrix in (11) cannot have two zero entries while the matrix in (10) has two zero entries.

We next prove that the assumptions of Theorem 1 hold for this example. Setting

$$k = \begin{bmatrix} 1 \\ 1/2 \\ 5/4 \end{bmatrix}, b = -1/4$$

we get

$$R = \left( diag(k) A + bkk^T A \right)^S = -\frac{1}{768} \begin{bmatrix} 600 & 66 & -77 \\ 66 & 108 & 197 \\ -77 & 197 & 470 \end{bmatrix}$$

It can be easily verified (using Sylvester's criterion) that $R$ is a negative definite matrix. It follows that all assumptions of Theorem 1 hold with $\kappa = \mu = 1$, $g = \delta = 0$, $p = q = \beta = 0$ and $F = \mathbb{R}^3_+$. Indeed, due to the fact that $R$ is a negative definite matrix the set $C$ defined by (9) coincides with the singleton $\{1_3\}$.

Since $bk^T 1_3 = -\frac{11}{16} > -1$, $b = -1/4 < 0$, $p = q = 0$, $k \in \operatorname{int}\left(\mathbb{R}^3_+\right)$ and since the Lotka-Volterra system (5) with $A$ given by (10) is permanent, all assumptions of Corollary 1 and Corollary 3 hold (with arbitrary $c^* > 0$). Therefore, $x^* = 1_3$ is globally attractive for (5) with $A$ given by (10). $\triangleleft$

We next discuss the relation of Theorem 1 and Corollaries 1, 2 and 3 with existing general results in the literature. More specifically, we show that the sufficient conditions for global attractivity that are provided by Corollaries 1, 2 and 3 are less demanding than the corresponding conditions of the existing general results in the literature.

I. Relation with Volterra-Lyapunov conditions

If there exists $h \in \operatorname{int}\left(\mathbb{R}^n_+\right)$ such that $\left(diag(h)A\right)^S$ is negative definite then the assumptions of Theorem 1 hold with $\kappa = \mu = 1$, $g = \delta = 0$, $k = h$, $b = 0$, $p = q = \beta = 0$ and $F = \mathbb{R}^3_+$. Indeed, due to the fact that $R$ is a negative definite matrix the set $C$ defined by (9) coincides with the singleton $\{1_n\}$. Therefore, all assumptions of Corollary 1 and Corollary 2 hold (with arbitrary $c^* > 0$).

II. Relation with the eigenvector conditions in [1]

In order to study the relation of Theorem 1 with the eigenvector conditions in [1], we need the following result

**Lemma 1 (see Definition 2.2.4 on page 14 and Theorem 2.3.3 on page 18 in [14]):** *Let two symmetric matrices* $Z_i \in \mathbb{R}^{n \times n}$, $i = 1,2$ *be such that* $\min\left(x^T Z_1 x, x^T Z_2 x\right) < 0$ *for every* $x \in \mathbb{R}^n \setminus \{0\}$. *Then there exist constants* $\tau_1, \tau_2 \geq 0$ *such that the matrix* $\left(\tau_1 Z_1 + \tau_2 Z_2\right)$ *is negative definite.*

If there exists $\alpha = \left(\alpha_1, ..., \alpha_n\right)^T \in \mathbb{R}^n$ and $\lambda < 0$ with $\alpha_i \neq 0$ for $i = 1, ..., n$ and $\alpha^T A = \lambda \alpha^T$ for which $y^T diag(\alpha) A y < 0$ for all $y \in \mathbb{R}^n$ with $y \neq 0$ and $\alpha^T y = 0$, then we can define the symmetric matrices $Z_1 = \left(diag(\alpha)A\right)^S$ and $Z_2 = -\alpha\alpha^T$. Clearly, $x^T Z_2 x = -\left(\alpha^T x\right)^2$ for all $x \in \mathbb{R}^n$. Therefore, we understand that for every $x \in \mathbb{R}^n \setminus \{0\}$ we either have $x^T Z_2 x < 0$ (when $\alpha^T x \neq 0$) or $x^T Z_2 x = 0$ (when $\alpha^T x = 0$) and $x^T Z_1 x = x^T \left(diag(\alpha)A\right)^S x = x^T diag(\alpha) A x < 0$. Therefore, $\min\left(x^T Z_1 x, x^T Z_2 x\right) < 0$ for every $x \in \mathbb{R}^n \setminus \{0\}$. Thus, by virtue of Lemma 1, there exist constants $\tau_1, \tau_2 \geq 0$ such that the matrix $\tau_1 Z_1 + \tau_2 Z_2 = \tau_1 \left(diag(\alpha)A\right)^S - \tau_2 \alpha\alpha^T$ is negative definite. We notice that we must have $\tau_1 > 0$, because if $\tau_1 = 0$ then the matrix $-\tau_2 \alpha\alpha^T$ would be negative definite; a contradiction since $-\tau_2 x^T \alpha\alpha^T x = -\tau_2 \left(\alpha^T x\right)^2 = 0$ for every $x \in \mathbb{R}^n \setminus \{0\}$ with $\alpha^T x = 0$. Since $\tau_1 > 0$, it follows that the matrix $\left(diag(\alpha)A - \frac{\tau_2}{\tau_1}\alpha\alpha^T\right)^S$ is negative definite. Since $\alpha^T = \frac{1}{\lambda}\alpha^T A$, it follows that the matrix $\left(diag(\alpha)A - \frac{\tau_2}{\lambda\tau_1}\alpha\alpha^T A\right)^S$ is negative definite. Setting

$$k = \alpha, b = -\frac{\tau_2}{\lambda\tau_1}$$

we get that $R = \left(diag(k)A + bkk^T A\right)^S$ is a negative definite matrix. Hence, the assumptions of Theorem 1 hold with $\kappa = \mu = 1$, $g = \delta = 0$, $k = \alpha$, $b = -\frac{\tau_2}{\lambda\tau_1}$, $p = q = \beta = 0$ and $F = \mathbb{R}^3_+$. Indeed, due to the fact that $R$ is a negative definite matrix the set $C$ defined by (9) coincides with the singleton $\{1_n\}$. Therefore, if the eigenvector conditions in [1] hold then all assumptions of Corollary 1 hold automatically.

It should be noted that the eigenvector conditions in [1] require permanence of the Lotka-Volterra system (5) while on the other hand Corollary 1 requires the less demanding property that for every $x_0 \in \operatorname{int}\left(\mathbb{R}^n_+\right)$ the solution $x(t)$ of (5) with initial condition $x(0) = x_0$ is bounded and satisfies $\liminf_{t\to+\infty}\left(x_i(t)\right) > 0$ for $i = 1,...,n$. No uniform bounds are needed for Corollary 1.

III. Relation with the results in [21]

The results in [21] deal with the case (5) with $n = 3$ where the matrix $A$ has the special form (11) with (12) and (13). It was shown in [21] that systems of this form include systems that fail to satisfy the Volterra-Lyapunov stability criteria and the eigenvector conditions in [1].

Using (11), (12), (13) it is straightforward to verify that (7) holds with

$$p = -\begin{bmatrix} y_1^* \\ y_2^* \\ y_3^* \end{bmatrix} \quad , \quad \beta = \begin{bmatrix} \frac{y_1^*}{\lambda_1} \\ \frac{y_2^*}{\lambda_2} \\ \frac{y_3^*}{\lambda_3} \end{bmatrix} \quad , \quad g = 0$$

Moreover, (6) and (8) hold with

$$\kappa = \lambda_2 > 0 \quad , \quad \mu = 1/\lambda_2 > 0$$

Selecting $k = q = 0$, $b = \delta = 0$, we notice that the matrices $R = \left(\left(diag\left(k\right) + bkk^T\right)A\right)^S$ and $Q = \left(diag\left(q\right)A - \delta pp^T\right)^S$ are negative semidefinite. Selecting $F = \mathbb{R}_+^3$, we next find the largest invariant set $J$ in the set

$$C := \left\{ x = \left(x_1, x_2, x_3\right)^T \in \mathbb{R}_+^3 : y_1^*\left(x_1 - 1\right) + y_2^*\left(x_2 - 1\right) + y_3^*\left(x_3 - 1\right) = 0 \right\} \tag{14}$$

We have (using (5) with $n = 3$, (11), (12), (13) and (14)) for every solution that evolves in $C$:

$$\begin{aligned} \dot{x}_1 &= -\frac{x_1}{\delta} y_3^*\left(x_3 - 1\right) \\ \dot{x}_2 &= 0 \\ \dot{x}_3 &= \frac{x_3}{\delta} y_1^*\left(x_1 - 1\right) \end{aligned} \tag{15}$$

Since any solution $x(t)$ in $C$ satisfies $\frac{d}{dt}\left(y_1^*\left(x_1(t) - 1\right) + y_2^*\left(x_2(t) - 1\right) + y_3^*\left(x_3(t) - 1\right)\right) \equiv 0$, using (15) we get $x_1(t) \equiv x_3(t)$. Finally, using the facts $\frac{d}{dt}\left(x_3(t) - x_1(t)\right) \equiv 0$, $x_1(t) \equiv x_3(t)$ we get $x_1(t)\left(x_1(t) - 1\right) \equiv 0$. Therefore, the largest invariant set $J$ in the set $C$ is the set

$$J = \left\{ \begin{bmatrix} 0 \\ \dfrac{1 - \lambda_2}{y_2^*} \\ 0 \end{bmatrix}, 1_3 \right\} \tag{16}$$

It follows from (16) that $\{1_3\} = \operatorname{int}\left(\mathbb{R}_+^3\right) \cap J$ and we conclude that all assumptions of Theorem 1 hold.

In order to prove global attractivity of $x^* = 1_3$ we can use Corollary 2. Indeed, $b = 0$, $(-p) \in \mathbb{R}_+^3$ and $q_i + k_i > c^* p_i$ for $i = 1, 2, 3$ and arbitrary $c^* > 0$. Moreover, we notice that the matrix $A$ given by (11), (12), (13) is a $B-matrix$ (this follows from Exercise 15.2.6 on page 188 and Theorem 15.2.11 on page 189 in [4] and the fact that $\det(-A) = \frac{y_1^* y_2^* y_3^*}{\lambda_2 \delta^2} > 0$). Hence, Theorem 15.2.1 on page 185 in [4] implies that the Lotka-Volterra system (5) with $n = 3$ and $A$ given by (11), (12), (13) is uniformly bounded. Therefore, if the stability conditions in [21] hold then all assumptions of Corollary 2 hold automatically.

We end this section by providing a Lotka-Volterra system for which Theorem 1 fails.

**Example 2:** System (5) with $n=3$ and

$$A=-\begin{bmatrix} 5 & 10 & 2 \\ 4 & 7 & 11 \\ 10 & 2 & 8 \end{bmatrix} \tag{17}$$

was studied in [22]. It was shown in [22] that it is impossible to satisfy the Volterra-Lyapunov conditions or the eigenvector conditions in [1] for system (5) with $A$ given by (17). It has been conjectured in [22] that $x^*=1_3$ is globally attractive for (5) with $A$ given by (17).

We have checked that it is not possible to write the matrix $A$ given by (17) in the form (11).

We have performed extensive numerical investigations and we have not been able to verify that system (5) with $A$ given by (17) satisfies the assumptions of Theorem 1. More specifically, straightforward calculations show that there are no constants $\kappa,\mu>0$, $g,\delta\geq 0$ and vectors $p\in\mathbb{R}^n$, $\beta\in\mathbb{R}_+^n$ other than $p=\beta=0$ for which (6), (7), (8) hold. Furthermore, our numerical investigations were not able to give a vector $k\in\mathbb{R}^n$ other than zero and a constant $b\in\mathbb{R}$ for which the matrix $R=\left(\left(diag\left(k\right)+bkk^T\right)A\right)^S$ is negative semidefinite. $\triangleleft$

## 4. Proof of Theorem 1

In this section we provide the proof of Theorem 1. For the proof of Theorem 1 we use LaSalle's theorem, which is stated next for reader's convenience.

**LaSalle's Theorem (see [9]):** *Let $\Omega\subseteq D$ be a compact and positively invariant set for the dynamical system $\dot{x}=f(x), x\in D$, where $D\subseteq\mathbb{R}^n$ is an open set and $f:D\to\mathbb{R}^n$ is a locally Lipschitz mapping. Let $U:D\to\mathbb{R}$ be a continuously differentiable function that satisfies $\dot{U}(x)\leq 0$ for all $x\in\Omega$. Then $\omega\left(x\right)\subseteq M$ for all $x\in\Omega$, where $M\subseteq\Omega$ is the largest invariant set in the set* $\left\{\xi\in\Omega:\dot{U}(\xi)=0\right\}$.

**Proof of Theorem 1:** We start by defining some useful functions and we calculate their time derivatives.

If $b\neq 0$ we define the function $V:\text{int}\left(\mathbb{R}_+^n\right)\to\mathbb{R}$ by means of the formula for all $x=\left(x_1,...,x_n\right)^T\in\text{int}\left(\mathbb{R}_+^n\right)$:

$$V(x):=\left(k^T\left(x-1_n\right)-\frac{1}{b}\right)\prod_{i=1}^{n}x_i^{bk_i}+\frac{1}{b} \tag{18}$$

Using the facts that $k^T diag(x) = x^T diag(k)$ and $k^T = 1_n^T diag(k)$, as well as (5) and definitions (18) and $R = \left(\left(diag(k) + bkk^T\right)A\right)^S$, we obtain for all $x \in \operatorname{int}\left(\mathbb{R}_+^n\right)$ when $b \neq 0$:

$$\begin{aligned}
\dot{V}(x) &= k^T diag(x) A(x-1_n)\prod_{i=1}^{n} x_i^{bk_i} + b\left(k^T(x-1_n) - \frac{1}{b}\right)k^T A(x-1_n)\prod_{i=1}^{n} x_i^{bk_i} \\
&= \left(x^T diag(k) A(x-1_n) + b\left(k^T(x-1_n)\right)k^T A(x-1_n) - k^T A(x-1_n)\right)\prod_{i=1}^{n} x_i^{bk_i} \\
&= \left(x^T diag(k) A(x-1_n) + b(x-1_n)^T kk^T A(x-1_n) - 1_n^T diag(k) A(x-1_n)\right)\prod_{i=1}^{n} x_i^{bk_i} \\
&= (x-1_n)^T \left(diag(k) + bkk^T\right) A(x-1_n)\prod_{i=1}^{n} x_i^{bk_i} \\
&= (x-1_n)^T \left(\left(diag(k) + bkk^T\right)A\right)^S (x-1_n)\prod_{i=1}^{n} x_i^{bk_i} = (x-1_n)^T R(x-1_n)\prod_{i=1}^{n} x_i^{bk_i}
\end{aligned} \tag{19}$$

If $b = 0$ we define the function $V : \operatorname{int}\left(\mathbb{R}_+^n\right) \to \mathbb{R}$ by means of the formula for all $x = (x_1,...,x_n)^T \in \operatorname{int}\left(\mathbb{R}_+^n\right)$:

$$V(x) := k^T(x-1_n) - \sum_{i=1}^{n} k_i \ln(x_i) \tag{20}$$

Using the facts that $k^T diag(x) = x^T diag(k)$ and $k^T = 1_n^T diag(k)$, as well as (5) and definitions (20) and $R = \left(\left(diag(k) + bkk^T\right)A\right)^S$, we obtain for all $x \in \operatorname{int}\left(\mathbb{R}_+^n\right)$ when $b = 0$:

$$\begin{aligned}
\dot{V}(x) &= k^T diag(x) A(x-1_n) - k^T A(x-1_n) \\
&= x^T diag(k) A(x-1_n) - 1_n^T diag(k) A(x-1_n) \\
&= (x-1_n)^T diag(k) A(x-1_n) \\
&= (x-1_n)^T \left(\left(diag(k) + bkk^T\right)A\right)^S (x-1_n) \\
&= (x-1_n)^T R(x-1_n)\prod_{i=1}^{n} x_i^{bk_i}
\end{aligned} \tag{21}$$

We conclude from (19) and (21) that for every $b \in \mathbb{R}$, the function $V : \operatorname{int}\left(\mathbb{R}_+^n\right) \to \mathbb{R}$ defined by (18) when $b \neq 0$ and (20) when $b = 0$ satisfies the following equation for all $x \in \operatorname{int}\left(\mathbb{R}_+^n\right)$:

$$\dot{V}(x) = (x-1_n)^T R(x-1_n)\prod_{i=1}^{n} x_i^{bk_i} \tag{22}$$

Next, we define for all $x \in \operatorname{int}\left(\mathbb{R}_+^n\right)$:

$$S = p^T\left(x - 1_n\right) + \kappa \tag{23}$$

Using the fact $p^T diag(x) = x^T diag(p)$, $1_n^T diag(p) = p^T$, (5), (7), (8) and (23) we get for all $x \in \operatorname{int}\left(\mathbb{R}_+^n\right)$:

$$\begin{aligned}
\dot{S} &= p^T diag(x) A(x - 1_n) = x^T diag(p) A(x - 1_n) \\
&= (x - 1_n)^T diag(p) A(x - 1_n) + 1_n^T diag(p) A(x - 1_n) \\
&= (x - 1_n)^T diag(p) A(x - 1_n) + p^T A(x - 1_n) \\
&= (x - 1_n)^T \left(diag(p) A + (\beta - gp) p^T\right)(x - 1_n) - (x - 1_n)^T (\beta - gp) p^T (x - 1_n) + p^T A(x - 1_n) \\
&= (x - 1_n)^T \left(diag(p) A + (\beta - gp) p^T\right)^S (x - 1_n) - (x - 1_n)^T (\beta - gp) p^T (x - 1_n) + p^T A(x - 1_n) \\
&= -(x - 1_n)^T (\beta - gp) p^T (x - 1_n) + p^T A(x - 1_n) \\
&= -(x - 1_n)^T \beta p^T (x - 1_n) + g(x - 1_n)^T p p^T (x - 1_n) - (\mu - g\kappa) p^T (x - 1_n) - \kappa\beta^T (x - 1_n) \\
&= -(S - \kappa)\beta^T (x - 1_n) + g(S - \kappa)^2 - (\mu - g\kappa)(S - \kappa) - \kappa\beta^T (x - 1_n)
\end{aligned} \tag{24}$$

It follows from (24) and (23) that the following equation holds for all $x \in \operatorname{int}\left(\mathbb{R}_+^n\right)$:

$$\dot{S} = -S\beta^T x + gS^2 - \left(\mu - \beta^T 1_n + g\kappa\right)S + \kappa\mu \tag{25}$$

We notice that since $\beta \in \mathbb{R}_+^n$, it follows that $\beta^T x \geq 0$ for all $x \in \operatorname{int}\left(\mathbb{R}_+^n\right)$. Therefore, (23), (25), inequality (6) and the fact that $g \geq 0$ guarantee the following implication for all $x \in \operatorname{int}\left(\mathbb{R}_+^n\right)$

$$S = p^T x + \kappa \leq \varepsilon \Rightarrow \dot{S} \geq \kappa\mu - \left(\beta^T x + \mu - \beta^T 1_n + g\kappa\right)\varepsilon \tag{26}$$

Moreover, it follows from (24) and (23) that the following equation holds for all $x \in \operatorname{int}\left(\mathbb{R}_+^n\right)$:

$$\dot{S} = -S(\beta - gp)^T (x - 1_n) - \mu(S - \kappa) \tag{27}$$

Therefore, we obtain from (23), (27) and (7) for all $x \in \operatorname{int}\left(\mathbb{R}_+^n\right)$ with $p^T(x - 1_n) + \kappa > 0$:

$$\frac{d}{dt}\left(S-\kappa-\kappa\ln\left(\frac{S}{\kappa}\right)\right)=\frac{S-\kappa}{S}\dot{S}$$
$$=-\left(x-1_n\right)^T\left(\beta-gp\right)p^T\left(x-1_n\right)-\frac{\mu}{p^T\left(x-1_n\right)+\kappa}\left(p^T\left(x-1_n\right)\right)^2$$
$$=-\left(x-1_n\right)^T\left(diag\left(p\right)A+\left(\beta-gp\right)p^T\right)\left(x-1_n\right) \tag{28}$$
$$-\frac{\mu}{p^T\left(x-1_n\right)+\kappa}\left(p^T\left(x-1_n\right)\right)^2+\left(x-1_n\right)^T diag\left(p\right)A\left(x-1_n\right)$$
$$=-\frac{\mu}{p^T\left(x-1_n\right)+\kappa}\left(p^T\left(x-1_n\right)\right)^2+\left(x-1_n\right)^T diag\left(p\right)A\left(x-1_n\right)$$

Let $c>0$ be an arbitrary positive constant. We define the function $W_c$ for all $x\in \operatorname{int}\left(\mathbb{R}^n_+\right)$ with $p^T\left(x-1_n\right)+\kappa>0$:

$$\begin{aligned}W_c(x)&=q^T\left(x-1_n\right)-\sum_{i=1}^{n}\left(q_i-cp_i\right)\ln\left(x_i\right)-c\kappa\ln\left(\frac{p^T\left(x-1_n\right)+\kappa}{\kappa}\right)\\&=(q-cp)^T\left(x-1_n\right)-\sum_{i=1}^{n}\left(q_i-cp_i\right)\ln\left(x_i\right)+c\left(S-\kappa-\kappa\ln\left(\frac{S}{\kappa}\right)\right)\end{aligned} \tag{29}$$

Using (28), (29), (5), the facts that $\left(q-cp\right)^T diag\left(x\right)=x^T diag\left(q-cp\right)$ and $(q-cp)^T=1_n^T diag\left(q-cp\right)$ as well as the definition $Q=\left(diag\left(q\right)A-\delta pp^T\right)^S$, we get for all $x\in \operatorname{int}\left(\mathbb{R}^n_+\right)$ with $p^T\left(x-1_n\right)+\kappa>0$:

$$\dot{W}_c(x)=(q-cp)^T diag(x)A\left(x-1_n\right)-(q-cp)^T A\left(x-1_n\right)+c\frac{d}{dt}\left(S-\kappa-\kappa\ln\left(\frac{S}{\kappa}\right)\right)$$
$$=x^T diag(q-cp)A\left(x-1_n\right)-1_n^T diag(q-cp)A\left(x-1_n\right)$$
$$-\frac{c\mu}{p^T\left(x-1_n\right)+\kappa}\left(p^T\left(x-1_n\right)\right)^2+c\left(x-1_n\right)^T diag\left(p\right)A\left(x-1_n\right)$$
$$=\left(x-1_n\right)^T diag\left(q\right)A\left(x-1_n\right)-\frac{c\mu}{p^T\left(x-1_n\right)+\kappa}\left(p^T\left(x-1_n\right)\right)^2 \tag{30}$$
$$=\left(x-1_n\right)^T\left(diag\left(q\right)A-\delta pp^T\right)\left(x-1_n\right)+\delta\left(p^T\left(x-1_n\right)\right)^2-\frac{c\mu}{p^T\left(x-1_n\right)+\kappa}\left(p^T\left(x-1_n\right)\right)^2$$
$$=\left(x-1_n\right)^T Q\left(x-1_n\right)-\left(\frac{c\mu}{p^T\left(x-1_n\right)+\kappa}-\delta\right)\left(p^T\left(x-1_n\right)\right)^2$$

Finally, we define the function $W_c$ for all $x\in \operatorname{int}\left(\mathbb{R}^n_+\right)$ with $p^T\left(x-1_n\right)+\kappa>0$:

$$U_c(x)=V(x)+W_c(x) \tag{31}$$

Combining (22), (31) and (30) we get for all $x \in \operatorname{int}\left(\mathbb{R}_+^n\right)$ with $p^T\left(x - 1_n\right) + \kappa > 0$:

$$\dot{U}_c(x) = \left(x - 1_n\right)^T R\left(x - 1_n\right) \prod_{i=1}^{n} x_i^{bk_i} + \left(x - 1_n\right)^T Q\left(x - 1_n\right) - \left(\frac{c\mu}{p^T\left(x - 1_n\right) + \kappa} - \delta\right)\left(p^T\left(x - 1_n\right)\right)^2 \quad (32)$$

Let arbitrary $x_0 \in \operatorname{int}\left(\mathbb{R}_+^n\right)$, for which the solution $x(t)$ of (5) with initial condition $x(0) = x_0$ is bounded, be given.

Let $\varepsilon > 0$ sufficiently small be such that

$$\rho := \kappa\mu - \left(\sup\left\{ \beta^T x(t) : t \geq 0 \right\} + \mu - \beta^T 1_n + g\kappa\right)\varepsilon > 0 \quad (33)$$

We notice that there exists $t_0 > 0$ such that $S(t) = p^T\left(x(t) - 1_n\right) + \kappa \geq \varepsilon$ for all $t \geq t_0$. Indeed, this fact is a direct consequence of (26) and (33), which guarantee the following implication:

$$\text{If } S(t) = p^T x(t) + \kappa \leq \varepsilon \text{ then } \dot{S}(t) \geq \rho > 0 \quad (34)$$

The assumptions of Theorem 1 for the set $F \subseteq \mathbb{R}_+^n$ allow us to conclude that there exists $\tau > 0$ such that $x(t) \in F$ for all $t \geq \tau$. Without loss of generality, we can assume that $S(t) = p^T\left(x(t) - 1_n\right) + \kappa \geq \varepsilon$ for all $t \geq \tau$.

For the proof of (a) we assume that the bounded solution $x(t)$ of (5) with initial condition $x(0) = x_0$ satisfies $\liminf_{t \to +\infty}\left(x_i(t)\right) > 0$ for $i = 1,...,n$. In this case we pick $c > 0$ so that

$$c \geq \frac{\delta + 1}{\mu}\left(\kappa + \sup\left\{ p^T\left(x(t) - 1_n\right) : t \geq 0 \right\}\right) \quad (35)$$

Equations (32), (35) and the facts that $x(t) \in \operatorname{int}\left(\mathbb{R}_+^n\right)$, $p^T\left(x(t) - 1_n\right) + \kappa \geq \varepsilon$ for all $t \geq \tau$ as well as and the fact that the matrices $R = \left(\left(diag\left(k\right) + bkk^T\right)A\right)^S$ and $Q = \left(diag\left(q\right)A - \delta pp^T\right)^S$ are negative semidefinite, allow us to conclude that the following inequality holds for all $t \geq \tau$:

$$\dot{U}_c(x(t)) \leq 0 \quad (36)$$

For the proof of (b) and (c) we show that the bounded solution $x(t)$ of (5) with initial condition $x(0) = x_0$ satisfies $\liminf_{t \to +\infty}\left(x_i(t)\right) > 0$ for $i = 1,...,n$. In these cases, we pick $c > 0$ so that

$$c \geq c^* + \frac{\delta + 1}{\mu}\left(\kappa + \sup\left\{ p^T\left(x(t) - 1_n\right) : t \geq 0 \right\}\right) \quad (37)$$

where $c^* > 0$ is the constant involved in case (b) or case (c). Again, equations (32), (37) and the facts that $x(t) \in \text{int}\left(\mathbb{R}_+^n\right)$, $p^T\left(x(t)-1_n\right)+\kappa \geq \varepsilon$ for all $t \geq \tau$ as well as and the fact that the matrices $R = \left(\left(diag\left(k\right)+bkk^T\right)A\right)^S$ and $Q = \left(diag\left(q\right)A - \delta pp^T\right)^S$ are negative semidefinite, allow us to conclude that inequality (36) holds for all $t \geq \tau$.

Claim 1: If $b = 0$, $(-p) \in \mathbb{R}_+^n$ and there exists $c^* > 0$ such that $q_i + k_i > c^* p_i$ for $i = 1,...,n$ then the bounded solution $x(t)$ of (5) with initial condition $x(0) = x_0$ satisfies $\liminf_{t\to+\infty}\left(x_i(t)\right) > 0$ for $i = 1,...,n$.

Proof of Claim 1: In this case we use (20) and (29) to obtain the following equation for all $x \in \text{int}\left(\mathbb{R}_+^n\right)$ with $p^T\left(x-1_n\right)+\kappa > 0$:

$$U_c(x) = \sum_{i=1}^{n}\left(q_i + k_i - cp_i\right) f\left(x_i\right) + c\kappa f\left(\frac{p^T\left(x-1_n\right)+\kappa}{\kappa}\right) \tag{38}$$

where $f\left(s\right) = s - 1 - \ln(s)$ for $s > 0$. The function $f : \left(0, +\infty\right) \to \mathbb{R}_+$ is a non-negative function with global minimum at $s = 1$ that satisfies $\lim_{s\to+\infty}\left(f(s)\right) = \lim_{s\to0^+}\left(f(s)\right) = +\infty$. Therefore, by showing that $f\left(x_i(t)\right)$ is bounded for $t \geq \tau$ and $i = 1,...,n$, we can conclude that $\liminf_{t\to+\infty}\left(x_i(t)\right) > 0$ for $i = 1,...,n$.

Since $c > c^*$ (recall (37)) and $p_i \leq 0$, $q_i + k_i > c^* p_i$ for $i = 1,...,n$ we obtain $q_i + k_i - cp_i > 0$ for $i = 1,...,n$. Inequality (36) implies that the following estimate holds for $t \geq \tau$:

$$U_c(x(t)) \leq U_c(x(\tau)) \tag{39}$$

The fact that the function $f : \left(0, +\infty\right) \to \mathbb{R}_+$ is a non-negative function allows us to combine (38) and (39) to get for $t \geq \tau$ and $i = 1,...,n$:

$$f\left(x_i(t)\right) \leq \left(q_i + k_i - cp_i\right)^{-1} U_c(x(\tau)) \tag{40}$$

Thus $f\left(x_i(t)\right)$ is bounded for $t \geq \tau$ and $i = 1,...,n$. The proof of the Claim 1 is complete. $\triangleleft$

Claim 2: If $k \in \mathbb{R}_+^n$, $(-p) \in \mathbb{R}_+^n$, $b < 0$ with $bk^T 1_n > -1$ and there exists $c^* > 0$ such that $q_i + k_i > c^* p_i$ and $q_i \geq c^* p_i$ for $i = 1,...,n$ then the bounded solution $x(t)$ of (5) with initial condition $x(0) = x_0$ satisfies $\liminf_{t\to+\infty}\left(x_i(t)\right) > 0$ for $i = 1,...,n$.

Proof of Claim 2: In this case we use (18) and (29) to obtain the following equation for all $x \in \text{int}\left(\mathbb{R}_+^n\right)$ with $p^T\left(x-1_n\right)+\kappa > 0$:

$$U_c(x) = \left( k^T (x - 1_n) - \frac{1}{b} \right) \exp\left( b \sum_{i=1}^{n} k_i \ln(x_i) \right) + \frac{1}{b} + \sum_{i=1}^{n} (q_i - cp_i) f(x_i) + c\kappa f\left( \frac{p^T (x - 1_n) + \kappa}{\kappa} \right) \quad (41)$$

where $f(s) = s - 1 - \ln(s)$ for $s > 0$.

Since $c > c^*$ (recall (37)) and $p_i \le 0$, $q_i \ge c^* p_i$, $q_i + k_i > c^* p_i$ for $i = 1,...,n$ we obtain $q_i - cp_i \ge 0$ and $q_i + k_i > cp_i$ for $i = 1,...,n$. Inequality (36) implies that estimate (39) holds for $t \ge \tau$.

Since $k \in \mathbb{R}_+^n$ and $x \in \mathrm{int}(\mathbb{R}_+^n)$ we get $k^T (x - 1_n) - \frac{1}{b} \ge -k^T 1_n - \frac{1}{b} = -\frac{1 + bk^T 1_n}{b}$ for all $x \in \mathrm{int}(\mathbb{R}_+^n)$. Therefore, we obtain for all $x \in \mathrm{int}(\mathbb{R}_+^n)$:

$$\left( k^T (x - 1_n) - \frac{1}{b} \right) \exp\left( b \sum_{i=1}^{n} k_i \ln(x_i) \right) \ge -\frac{1 + bk^T 1_n}{b} \exp\left( b \sum_{i=1}^{n} k_i \ln(x_i) \right) \quad (42)$$

Using the fact that the function $f : (0, +\infty) \to \mathbb{R}_+$ is a non-negative function and since $q_i - cp_i \ge 0$ for $i = 1,...,n$, we obtain from (41) and (42) and the fact that $b < 0$ and $bk^T 1_n + 1 > 0$ for all $x \in \mathrm{int}(\mathbb{R}_+^n)$ with $p^T (x - 1_n) + \kappa > 0$:

$$U_c(x) - \frac{1}{b} \ge -\frac{1 + bk^T 1_n}{b} \exp\left( b \sum_{i=1}^{n} k_i \ln(x_i) \right) > 0 \quad (43)$$

Since $b < 0$ and $bk^T 1_n + 1 > 0$ we get from (43) for all $x \in \mathrm{int}(\mathbb{R}_+^n)$ with $p^T (x - 1_n) + \kappa > 0$:

$$\frac{1}{b} \ln\left( \frac{1 - bU_c(x)}{1 + bk^T 1_n} \right) \le \sum_{i=1}^{n} k_i \ln(x_i) \quad (44)$$

Combining (44) and (39) and using the fact that $b < 0$ and $bk^T 1_n + 1 > 0$, we get for $t \ge \tau$:

$$\frac{1}{b} \ln\left( \frac{1 - bU_c(x(\tau))}{1 + bk^T 1_n} \right) \le \frac{1}{b} \ln\left( \frac{1 - bU_c(x(t))}{1 + bk^T 1_n} \right) \le \sum_{i=1}^{n} k_i \ln(x_i(t)) \quad (45)$$

Since $x(t)$ is bounded, there exists $G > 0$ such that $x_i(t) \le G$ for all $t \ge 0$ and $i = 1,...,n$. Therefore, since $k \in \mathbb{R}_+^n$ we get from (45) for $t \ge \tau$ and $i = 1,...,n$:

$$\frac{1}{b} \ln\left( \frac{1 - bU_c(x(\tau))}{1 + bk^T 1_n} \right) - \sum_{j \ne i} k_j \ln(G) \le k_i \ln(x_i(t)) \quad (46)$$

Using the fact that the function $f:(0,+\infty)\to\mathbb{R}_+$ is a non-negative function and since $q_i - cp_i \ge 0$ for $i=1,...,n$, we obtain from (41), (42), (39) and the fact that $b<0$ and $bk^T 1_n + 1 > 0$ for all $t \ge \tau$ and $i=1,...,n$:

$$(q_i - cp_i)\, f\left(x_i(t)\right) \le U_c(x(\tau)) - \frac{1}{b} \tag{47}$$

If $q_i - cp_i > 0$ for some $i=1,...,n$ then (47) implies that $f\left(x_i(t)\right)$ is bounded for $t \ge \tau$. Since $f:(0,+\infty)\to\mathbb{R}_+$ is a non-negative function with global minimum at $s=1$ that satisfies $\lim\limits_{s\to+\infty}\left(f(s)\right) = \lim\limits_{s\to 0^+}\left(f(s)\right) = +\infty$, we can conclude that $\liminf\limits_{t\to+\infty}\left(x_i(t)\right) > 0$.

If $q_i - cp_i = 0$ for some $i=1,...,n$ then the fact that $q_i + k_i > cp_i$ implies that $k_i > 0$. Therefore, (46) allows us to conclude that $\liminf\limits_{t\to+\infty}\left(x_i(t)\right) > 0$.

The proof of the Claim 2 is complete. $\triangleleft$

Thus, the bounded solution $x(t)$ of (5) with initial condition $x(0) = x_0$ satisfies $\liminf\limits_{t\to+\infty}\left(x_i(t)\right) > 0$ for $i=1,...,n$.

We next define the positive orbit of $x(\tau) \in \operatorname{int}\left(\mathbb{R}^n_+\right)$ by the formula:

$$\gamma^+ := \left\{\, x(t) : t \ge \tau \,\right\} \tag{48}$$

Since the solution $x(t)$ of (5) with initial condition $x(0) = x_0$ is bounded and satisfies $\liminf\limits_{t\to+\infty}\left(x_i(t)\right) > 0$ for $i=1,...,n$, it follows that the closure $\overline{\gamma^+}$ of the set $\gamma^+$ is a compact set that satisfies $\overline{\gamma^+} \subseteq \operatorname{int}\left(\mathbb{R}^n_+\right)$. Moreover, since $x(t) \in F$ and $p^T\left(x(t) - 1_n\right) + \kappa \ge \varepsilon$ for all $t \ge \tau$ and since $F \subseteq \mathbb{R}^n_+$ is closed, it follows that $\overline{\gamma^+} \subseteq \left\{\, \xi \in F : p^T\xi + \kappa \ge \varepsilon \,\right\}$. Finally, we notice that definition (48) (which implies that $\gamma^+$ is positively invariant) and the fact that the omega limit set $\omega\left(x(\tau)\right)$ of $x(\tau) \in \operatorname{int}\left(\mathbb{R}^n_+\right)$ is invariant, guarantee that $\overline{\gamma^+} = \gamma^+ \cup \omega\left(x(\tau)\right)$ is positively invariant.

Since (36) and definition (48) imply that $\dot{U}_c(\xi) \le 0$ for all $\xi \in \gamma^+$, continuity of $\dot{U}_c$ on $\left\{\, \xi \in \operatorname{int}\left(\mathbb{R}^n_+\right) : p^T\xi + \kappa > 0 \,\right\}$ and the facts that $\overline{\gamma^+} \subseteq \left\{\, \xi \in F : p^T\xi + \kappa \ge \varepsilon \,\right\}$ and $\overline{\gamma^+} \subseteq \operatorname{int}\left(\mathbb{R}^n_+\right)$ guarantee that $\dot{U}_c(\xi) \le 0$ for all $\xi \in \overline{\gamma^+}$.

Applying LaSalle's theorem for (5) with $D = \left\{\, \xi \in \operatorname{int}\left(\mathbb{R}^n_+\right) : p^T\xi + \kappa > 0 \,\right\}$ on the compact and positively invariant set $\Omega = \overline{\gamma^+} \subseteq \left\{\, \xi \in F \cap \operatorname{int}\left(\mathbb{R}^n_+\right) : p^T\xi + \kappa \ge \varepsilon \,\right\}$ allows us to conclude that $\omega(\xi) \subseteq M$ for all $\xi \in \overline{\gamma^+}$, where $M \subseteq \overline{\gamma^+}$ is the largest invariant set in $\left\{ \zeta \in \overline{\gamma^+} : \dot{U}_c(\zeta) = 0 \right\}$. Equation (32), inequality (35), the fact that $\sup\left\{ p^T \left(x(t) - 1_n\right) : t \geq 0 \right\} \geq \max\left\{ p^T \left(\xi - 1_n\right) : \xi \in \overline{\gamma^+} \right\}$ (a consequence of definition (48)) and the fact that the matrices $R = \left(\left(diag\left(k\right) + bkk^T\right)A\right)^S$ and $Q = \left(diag\left(q\right)A - \delta pp^T\right)^S$ are negative semidefinite imply that $M \subseteq \left\{ x \in \overline{\gamma^+} : \left(x - 1_n\right)^T R\left(x - 1_n\right) = \left(x - 1_n\right)^T Q\left(x - 1_n\right) = p^T\left(x - 1_n\right) = 0 \right\}$.

Since $\overline{\gamma^+} \subseteq \left\{ \xi \in F : p^T \xi + \kappa \geq \varepsilon \right\} \subseteq F$, it follows that $M \subseteq C$, where $C$ is the set defined in (9). Therefore, $M \subseteq J$ where $J \subseteq F$ is the largest invariant set in $C$. Since $M \subseteq \overline{\gamma^+} \subseteq \operatorname{int}\left(\mathbb{R}_+^n\right)$, it follows that $M \subseteq \operatorname{int}\left(\mathbb{R}_+^n\right) \cap J$. Since $\left\{ 1_n \right\} = \operatorname{int}\left(\mathbb{R}_+^n\right) \cap J$, we get that $M \subseteq \left\{ 1_n \right\}$.

The fact that $\omega(\xi) \subseteq M$ for all $\xi \in \overline{\gamma^+}$, the fact that $\omega(\xi) \neq \emptyset$ for all $\xi \in \overline{\gamma^+}$ (a consequence of the fact that $\overline{\gamma^+}$ is compact and positively invariant) and the fact that $M \subseteq \left\{ 1_n \right\}$ imply that $\omega(\xi) = \left\{ 1_n \right\}$ for all $\xi \in \overline{\gamma^+}$.

Since $x(\tau) \in \gamma^+ \subseteq \overline{\gamma^+}$ (recall definition (48)), it follows that $\omega\left(x(\tau)\right) = \left\{ 1_n \right\}$. The fact that $\omega\left(x_0\right) = \omega\left(x(0)\right) = \omega\left(x(\tau)\right)$ guarantees that $\omega\left(x_0\right) = \left\{ 1_n \right\}$. The proof is complete. $\lhd$

## 5. Conclusions

The present paper provided global attractivity results for the interior equilibrium point of a general Lotka-Volterra system with no restriction on the dimension of the system and with no special structure or properties of the interaction matrix. Theorem 1 contains as special cases all known general results. Moreover, global attractivity of the interior equilibrium point was shown for the three-dimensional Example 1, where none of the existing general results can be applied.

Example 2 as well as other results for competitive systems show that there is a lot of room for improvement. A possible direction that can be followed is to assume that the matrix $A$ is Hurwitz and use the recent results for $\omega-$limit sets in [8]. Another direction that can be followed is the use of the set $F \subseteq \mathbb{R}_+^n$ mentioned in Theorem 1. A possible use of the idea of ultimate contracting cells (see [5, 6]) for the construction of useful sets $F \subseteq \mathbb{R}_+^n$ may provide less demanding results that guarantee global attractivity of the interior equilibrium point.

## Acknowledgments

The author would like to thank Dr. Dionysis Theodosis for his help in Example 1 and Example 2 of the paper.